\documentclass[aos,seceqn,citesort,dvips]{arximspdf}
\usepackage{graphicx}

\doi{10.1214/10-AOS839}
\volume{39}
\issue{1}
\pubyear{2011}
\firstpage{278}
\lastpage{304}

\makeatletter

\newtheorem{theorem}{Theorem}
\newtheorem{prp}{Proposition}
\newtheorem{cor}{Corollary}

\newcommand{\varliminf}{\mathop{\underline{\lim}}}
\newcommand{\R}{\mathbb{R}}

\newcommand{\PP}{\mathbb{P}}
\newcommand{\E}{\mathbb{E}}

\newcommand{\cA}{\mathcal{A}}
\newcommand{\cF}{\mathcal{F}}

\newcommand{\eps}{\varepsilon}

\newcommand{\cS}{\mathcal{S}}
\newcommand{\cH}{\mathcal{H}}
\newcommand{\cN}{\mathcal{N}}
\newcommand{\cK}{\mathcal{K}}
\newcommand{\bbZ}{{\mathbb{Z}}}

\newcommand{\bbR}{\mathbb{R}}

\newcommand{\im}{\operatorname{im}}

\newcommand{\V}{\mathbb{V}}
\newcommand{\T}{\mathbb{T}}

\newcommand{\var}{\operatorname{Var}}

\newcommand{\ph}{\varphi}

\newcommand{\Zd}{\mathbb{Z}^d}
\newcommand{\Km}{\cK_m}
\newcommand{\Kmt}{\Km}
\newcommand{\Kt}{K}
\newcommand{\hyp}{\mathrm{H}}
\newcommand{\oB}{B}
\newcommand{\Lambdamax}{\overline{\Lambda}}
\newcommand{\Lambdamin}{\underline{\Lambda}}

\makeatother

\begin{document}
\begin{frontmatter}

\title{Detection of an anomalous cluster in a network\thanksref{T1}}
\runtitle{Detection of an anomalous cluster in a network}

\thankstext{T1}{Supported in part by NSF Grant DMS-06-03890 and ONR Grant
N00014-09-1-0258.}

\begin{aug}
\author[A]{\fnms{Ery} \snm{Arias-Castro}\corref{}\ead[label=e1]{eariasca@ucsd.edu}},
\author[B]{\fnms{Emmanuel J.} \snm{Cand\`es}\ead[label=e2]{candes@stanford.edu}} and
\author[C]{\fnms{Arnaud} \snm{Durand}\ead[label=e3]{arnaud.durand@math.u-psud.fr}}
\runauthor{E. Arias-Castro, E. J. Cand\`es and A. Durand}
\affiliation{University of California, San Diego, Stanford
University\break and Universit\'e Paris-Sud 11}
\address[A]{E. Arias-Castro\\
Department of Mathematics\\
University of California, San Diego\\
La Jolla, California 92093\\
USA\\
\printead{e1}}
\address[B]{E. J. Cand\`es\\
Departments of Mathematics\\
\quad and Statistics\\
Stanford University\\
Stanford, California 94305\\
USA\\
\printead{e2}}
\address[C]{A. Durand\\
Laboratoire de Math\'ematiques, UMR8628\\
Universit\'e Paris-Sud 11\\
Orsay, F-91405\\
France\\
\printead{e3}}
\end{aug}

\received{\smonth{1} \syear{2010}}
\revised{\smonth{6} \syear{2010}}

%
\begin{abstract}
We consider the problem of detecting whether or not, in a given sensor
network, there is a cluster of sensors which exhibit an ``unusual
behavior.'' Formally, suppose we are given a set of nodes and attach a
random variable to each node. We observe a realization of this process
and want to decide between the following two hypotheses: under the
null, the variables are i.i.d. standard normal; under the alternative,
there is a cluster of variables that are i.i.d. normal with positive
mean and unit variance, while the rest are i.i.d. standard normal.
%
We also address surveillance settings where each sensor in the network
collects information over time. The resulting model is similar, now
with a time series attached to each node. We again observe the process
over time and want to decide between the null, where all the variables
are i.i.d. standard normal, and the alternative, where there is an
emerging cluster of i.i.d. normal variables with positive mean and unit
variance. The growth models used to represent the emerging cluster are
quite general and, in particular, include cellular automata used in
modeling epidemics.
In both settings, we consider classes of clusters that are quite
general, for which we obtain a lower bound on their respective minimax
detection rate and show that some form of scan statistic, by far the
most popular method in practice, achieves that same rate to within a
logarithmic factor. Our results are not limited to the normal location
model, but generalize to any one-parameter exponential family when the
anomalous clusters are large enough.
\end{abstract}

%

\end{frontmatter}

\section{Introduction}

We discuss the problem of detecting whether or not, in a given network,
there is a cluster of nodes which exhibit an ``unusual behavior.''
Suppose that we are given a set of nodes with a random variable
attached to each node.
We observe a realization of this process and would like to tell whether
all the variables at the nodes have the same behavior, in the sense
that they are all sampled from a common distribution, or whether there
is a cluster of nodes at which the variables have a different distribution.

\subsection{A wide array of applications}

The task of detection in networks is critical for an increasing number
of applications, for example, in surveillance and environment
monitoring. We describe a few of these applications below.

\subsubsection*{Detection in sensor networks}
The advent of sensor networks
\cite{overviewsensor,1024422,zhao2004wireless} has multiplied the
amount of
data and the variety of applications where the task of detection is central.
Surveillance and environment monitoring are prime areas of application
for sensor networks.
Take, for example, the transport of hazardous materials.
Currently, some major traffic bottlenecks (e.g., airports, subways and
borders) use portal monitoring systems \cite{fitch2003tcr,1352095}.
Sensor networks offer a more flexible, decentralized alternative and
are considered for the detection of radioactive, biological or chemical
materials \cite{hills2001sd,radiation,cui2001nnh}.
Sensor networks are also extensively used in other target tracking
settings \cite{linesand,targetsensor}.

\subsubsection*{Detection in digital signals and images}
A digital camera may
be seen as a sensor network, with CCD or CMOS pixel sensors. As imaging
systems have been available for quite some time, the literature on
detection in images is quite extensive, spanning several decades,
particularly in satellite imagery
\cite{roadtracking,artificialnatural,shipdetection,firedetection},
computer vision \cite{eyedetection,objecttracking} and medical
imaging
\cite{medicalsurvey,breasttumor,braintumor,braams1995detection}.

\subsubsection*{Disease outbreak detection}
The presence of a biological or chemical material in a given
geographical region may also be detected indirectly through its impact
on human health.
In this context, early detection of the disease outbreak is crucial in
order to minimize the severity of the epidemic.
For that purpose, some specific information networks are used, with
surveillance systems now incorporating data from hospital emergency
visits, ambulance dispatch calls and pharmacy sales of over-the-counter
drugs \cite{rotz2004advances,heffernan2004ssp,wagner2001esv}.

\subsubsection*{Virus detection in a computer network}
Diseases affect computers as well, in the form of viruses and worms
spreading from host to host in a computer network~\cite{szor2005acv}.
Affected machines may exhibit slightly anomalous behavior (e.g., a~loss
of performance or violations of specific rules) which may be hard to
detect on an individual machine.

\subsubsection*{Detection from field measurements}
In \cite{waterquality}, the water quality in a network of streams in
Pennsylvania is assessed by field biologists performing a variety of
analyses at various locations along the streams; the objective is to
determine whether there are regions of low biological integrity based
on the collected data, and to identify these regions.
Other field measurements include census data and surveys involving
geographical location.

Detection is, of course, closely related to estimation (i.e., the
localization or
extraction of the anomalous cluster of nodes), but different. This
distinction is rarely
made clear, however. Indeed, reliable detection is possible at lower
signal-to-noise ratios than reliable estimation and it may be important
to detect the presence of signals from noisy data without being able to
estimate them.
For example, one could imagine developing a surveillance system
performing detection at relatively low energy/bandwidth costs, yet
efficient at low signal-to-noise ratios, and then switching to
estimation mode whenever the presence of a signal is detected. Another
example would be a low cost preliminary survey involving fewer field
measurements, with findings subsequently confirmed by a larger, more
expensive survey.

\subsection{Mathematical framework}

\subsubsection{Purely spatial model}
\label{sec:static}

We loosely model a network with a set of $m$ nodes, denoted by $\V_m$.
In our examples, we will either assume that $\V_m$ is embedded in a
Euclidean space or we will equip $\V_m$ with a graph structure. Our
analysis is in the setting of large networks, that is, $m \to\infty$.
To each node $v \in\V_m$, we attach a random variable $X_v$. The nodes
represent the sources of information (e.g., sensors) and the variables
represent the data they collect. In some settings, the data collected
by each unit is multidimensional, in which case $X_v$ is a random
vector. Our discussion readily generalizes to that setting.

The random variables are assumed to be independent. For concreteness,
we consider a normal location model, popular in signal and image
processing, to model the noise. Our analysis, however, generalizes to
any exponential family under some conditions on the sizes of the
anomalous clusters, such as Bernoulli models which arise in sensor
arrays where each sensor collects one bit (i.e., makes a binary
decision) or Poisson models which come up with count data, for
instance, arising in infectious disease surveillance systems
\cite{Kul}. The extension to exponential families is detailed in
Section \ref{sec:expo}.

\begin{figure}

\includegraphics{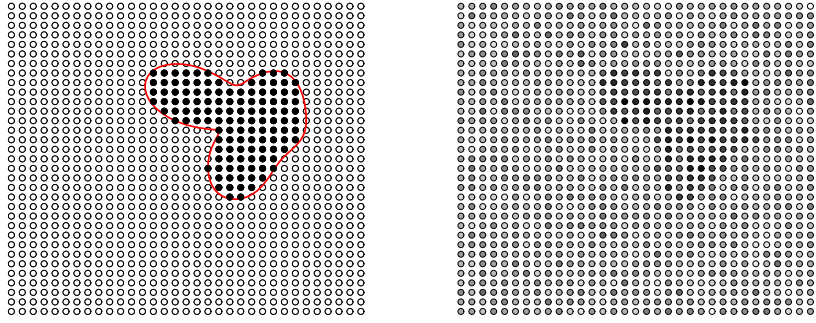}

\caption{Left: a thick cluster is defined as the nodes
within a closed curve, which is a mild deformation of a circle. Right:
corresponding noisy data.} \label{fig:thick}
\end{figure}

\begin{figure}[b]

\includegraphics{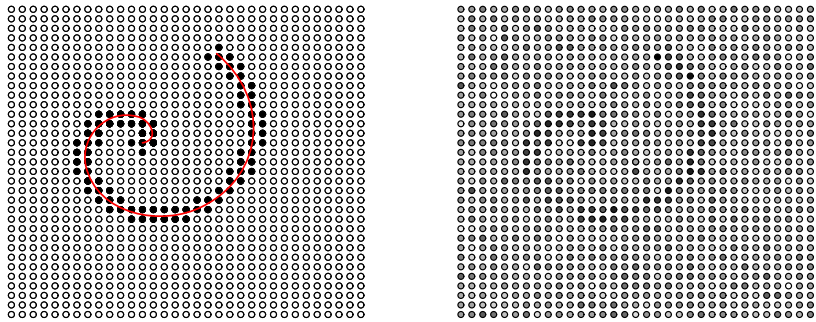}

\caption{Left: a thin cluster is defined as the nodes within
a band around a given curve. Right: corresponding noisy data.}
\label{fig:thin}
\end{figure}

\begin{figure}

\includegraphics{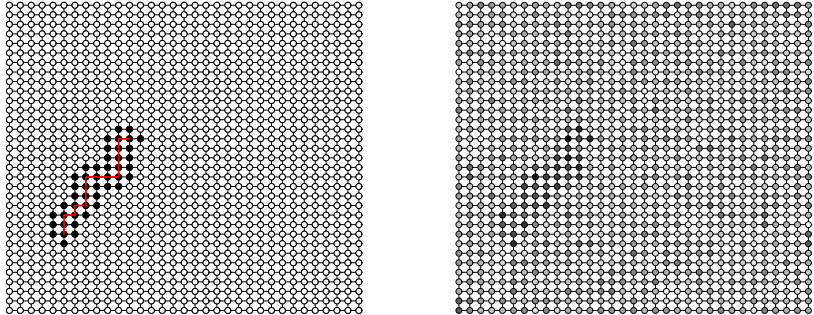}

\caption{Left: a band defined around a path. Right:
corresponding noisy data.} \label{fig:band}
\end{figure}

\begin{figure}[b]

\includegraphics{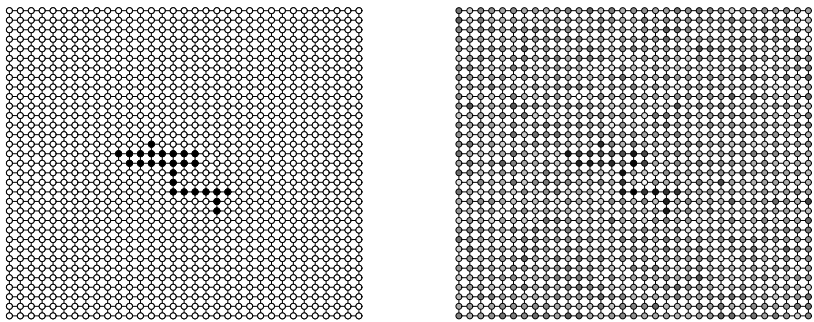}

\caption{Left: an arbitrary connected component. Right:
corresponding noisy data.} \label{fig:animal}
\end{figure}

%
%
%
The situation where no signal is present, that is, ``business as
usual,'' is modeled as
\[
\hyp_0^m \dvtx X_v \sim\cN(0,1)\qquad \forall v \in\V_m.
\]
Let $K$ be a cluster, which we define for now as a subset of nodes,
that is, \mbox{$K \subset\V_m$}. In fact, we will be interested in classes of
clusters that are either derived from a geometric shape, when $\V_m$ is
embedded in Euclidean space, or connected components, when $\V_m$ has a
graph structure.
The situation where the nodes in $K$ behave anomalously is modeled as
\[
\hyp_{1,K}^m \dvtx X_v \sim\cN(\mu_K,1)\qquad \forall v \in K;\qquad X_v \sim
\cN(0,1)\qquad \forall v \notin K,
\]
where $\mu_K > 0$.
We choose to decompose $\mu_K$ as $\mu_K = |K|^{-1/2} \Lambda_K$, where
$|K|$ denotes the number of nodes in $K$ and $\Lambda_K$ is the signal
strength. Indeed, with this normalization, for any cluster $K$,
%
%
\begin{equation} \label{eq:simple-risk}
\min_T \PP(T = 1|\hyp_0^m) + \PP(T = 0|\hyp_{1,K}^m) = 2 \PP\bigl(\cN
(0,1) >
\Lambda_K/2\bigr),
\end{equation}
where the minimum is over all tests for $\hyp_0^m$ versus $\hyp
_{1,K}^m$ and the lower bound is achieved by the likelihood ratio
(Neyman--Pearson) test. We define
\[
\Lambdamax_m = \max_{K \in\Km} \Lambda_K,\qquad \Lambdamin_m =
\min_{K
\in\Km} \Lambda_K.
\]
Figures \ref{fig:thick}--\ref{fig:animal} illustrate the setting for
various types of clusters.

%
%
%

%
%

Let $\Km$ be a class of clusters within $\V_m$ and define
\[
\hyp_1^m = \bigcup_{K \in\Km} \hyp_{1,K}^m.
\]
We are interested in testing $\hyp_0^m$ versus $\hyp_1^m$. In other
words, under the alternative, the cluster of anomalous nodes is only
known to belong to $\Km$.
We adopt a minimax point of view.
For a test $T$, we define its worst-case risk as
\[
\gamma_{\Km}(T) = \PP(T = 1|\hyp_0^m) + \max_{K \in\Km} \PP(T =
0|\hyp_{1,K}^m).
\]
The minimax risk for $\hyp_0^m$ versus $\hyp_1^m$ is defined as
\[
\gamma_{\Km} = \inf_T \gamma_{\Km}(T).
\]
We say that $\hyp_0^m$ and $\hyp_1^m$ are asymptotically inseparable
(in the minimax sense) if
\[
\varliminf_{m \to\infty} \gamma_{\Km} = 1,
\]
which is equivalent to saying that, as $m$ becomes large, no test can
perform substantially better than random guessing, without even looking
at the data.
A sequence of tests $(T_m)$ is said to asymptotically separate $\hyp
_0^m$ and $\hyp_1^m$ if
\[
\lim_{m \to\infty} \gamma_{\Km}(T_m) = 0,
\]
and $\hyp_0^m$ and $\hyp_1^m$ are said to be asymptotically separable
if there is such a sequence of tests.
For example, in view of (\ref{eq:simple-risk}), for any sequence of
clusters $K_m \subset\V_m$, $\hyp_0^m$ and $\hyp_{1, K_m}^m$ are
asymptotically inseparable if $\Lambda_{K_m} \to0$ and they are
asymptotically separable if $\Lambda_{K_m} \to\infty$.
For convenience, we assume that no cluster in the class $\Km$ is of
size comparable to that of the entire network, that is, $\max\{|K|\dvtx K
\in\Km\} = o(m)$. This simplifies the statement of our results and
detecting such clusters can easily be achieved using the test that
rejects for large values of $\sum_{v \in\V_m} X_v$.

%
%
%

%
%

The situation we just described is purely spatial and relevant in some
applications not involving time. Such situations are common in image
processing. In other applications, especially in surveillance, time is
an intrinsic part of the setting. In the following section, we modify
the model above to incorporate time.

\subsubsection{Spatio-temporal model}
\label{sec:dynamic}

Building on the framework introduced in the previous section, we assume
that each $X_v$ is now a (discrete) time series, $(X_v(t), t \in\T
_m)$, where $\T_m \subset[0, \infty)$ is finite with $|\T_m| \to
\infty
$; let $t_m = \max\{t \in\T_m\}$.
Let $\Kmt$ be a class of cluster sequences of the form $(K_t, t \in\T
_m)$ such that $K_t \subset\V_m$ for all $t \in\T_m$. For example,
assuming that $\V_m$ is embedded in a Euclidean space, with norm
denoted by $\| \cdot\|$, a space--time cylinder (e.g., one used in
disease outbreak detection \cite{diseaseoutbreak}) is a cluster
sequence $(K_t, t\in\mathbb{T}_m)$ of the form $K_t=\{v\in\mathbb
{V}_m\dvtx\|v-x_0\|\leq r_0\}$ if $t\geq t_0$, and $K_t=\varnothing$
otherwise, so that $t_0$ is the origin of the cluster in time and $x_0$
its center. Note that the radius remains constant here.
Another example is that of a space--time cone, of the form $K_t = \{v
\in\V_m\dvtx\|v - x_0\| \leq C(t-t_0)\}$ if $t\geq t_0,$ and
$K_t=\varnothing$ otherwise, so that $(x_0,t_0)$ is the origin of the
cluster in space--time.
The random variables $\{X_v(t)\dvtx v \in\V_m, t \in\T_m\}$ are assumed
to be independent.
This spatio-temporal setting is a special case of the purely spatial
setting with the set of nodes $\V_m \times\T_m$.
Understood as such, we are interested in testing $\hyp_0^m$ versus
$\hyp
_1^m$ as before.

\subsection{Structured multiple hypothesis testing}
\label{sec:multi-hypo}

Although the detection problem formulated above seems of great practical
relevance, the statistics literature is almost silent on the subject,
with the notable exception of the closely related topics of
change-point analysis \cite{brodsky1993nmc} and sequential
analysis \cite{MR799155}.
Indeed, the former is a special case of the spatial setting with the
one-dimensional lattice, while the latter is a special case of the
spatio-temporal setting where $\V_m$ has only one node.
In our context, these two settings are actually equivalent.

What is further puzzling is that a number
of publications addressing the task of detection in sensor networks
all assume overly simplistic models. For example, in
\cite{decentralized,energydriven,hierarchicalcensoring,optimaldistributed,cooperate,MR2450850},
the values at the sensors are
assumed to all have the same distribution under the null and the
alternative. That is, either all of the nodes are all right or
they are all anomalous---in our notation, $\Km= \{ \V_m \}$. First,
this is not a subtle statistical problem since, in such circumstances, it
suffices to apply the optimal likelihood ratio test.
Second, this assumption does not make sense in all of the applications
described above, where the event to be detected is expected to only
affect a small fraction of locations in the network.

In stark contrast, in all of the applications described earlier, the
set of alternatives is composite. Viewing each node as performing a
test of hypotheses, which is common in the literature on sensor
networks, our problem falls within the framework of multiple comparisons.
Multiple hypothesis testing is a rich and active line of research which is
receiving a considerable amount of attention within the statistical
community at the moment; see \cite{MR2387976} and references therein.
The vast majority of the papers assume that the tests are independent
of each other, which is clearly not the case here since, in general,
the class contains clusters that intersect. This is particularly true
in engineering applications, although this assumption is often
made \cite{MR2213561,MR926001,MR1951265}.

\subsection{The scan statistic}

We will focus on the test that rejects for large values of the
following version of the scan statistic:
%
%
\begin{equation}
\label{eq:scan}
\max_{K \in\Km} \frac{1}{\sqrt{|K|}} \sum_{v \in K} X_v.
\end{equation}
The chosen normalization is such that each term in the maximization is
standard normal under the null and allows us to compare clusters of
different sizes. It corresponds to the generalized likelihood ratio
test in our context if $\Lambda_K$ is independent of $K \in\Km$.
The scan statistic was originally proposed in the context of cluster
detection in point clouds \cite{glaz01}.
This is the method of matched filters which is ubiquitous in problems
of detection in a wide variety of fields, sometimes in the form of
deformable templates in the engineering literature
\cite{deformablereview,medicalsurvey} or their nonparametric equivalent,
active contours or snakes \cite{xu1998ssa}.
Note that the scan statistic is the prevalent method in disease
outbreak detection, with many variations
\cite{kulldorff2005stp,treebased,kulldorff2006ess,duczmal2006ess}.

As advocated in \cite{MGD}, we will not use the scan statistic directly
in most cases, but rather restrict the scanning to a subset of $\Km$.
More precisely, we will introduce, on subsets of nodes $K, L \subset\V_m$,
the metric
%
%
\begin{equation} \label{eq:delta}
\delta(K, L) = \sqrt{2} \biggl(1 - \frac{|K \cap L|}{\sqrt{|K|
|L|}}\biggr)^{1/2}
\end{equation}
and will restrict the scanning to an $\eps$-net of $\Km$ with respect
to $\delta$, that is, a subset $\{K_j\dvtx j \in J\} \subset\Km$ with the
property that for each $K \in\Km$, there is a $j \in J$ such that
$\delta(K, K_j) \leq\eps$. We will elaborate on this approach
in the \ref{supp}. When $J$ is minimal, we call the resulting statistic an $\eps
$-scan statistic. The approximation precision $\eps$ will be chosen
appropriately, depending on the situation.

We focus on $\eps$-scan statistics for two reasons. First, their
performance is easier to analyze than that of the scan statistic
itself; in fact, the main approach to analyzing the scan statistic, the
chaining method of Dudley \cite{MR0220340,MR2133757}, is via a properly
chosen $\eps$-scan statistic. Second, some of the classes we consider
are rather large and we believe that it would be computationally
impractical to scan through all of the clusters in the class;
furthermore, our results show that, from an asymptotic standpoint, no
substantial improvement would be gained by using the full scan statistic.

We also note that the tuning parameter $\eps$ may be dispensed with if
we scan over subsets of different sizes in a multiscale fashion and use
a scale-dependent threshold.


\subsection{Existing theoretical results}

The vast majority of the literature assumes that the set of nodes is
embedded in some Euclidean space, that is, $\V_m \subset\R^d$. This is
the case when the nodes represent spatial locations, such as in most
sensor networks. In this context, the cluster class $\Km$ is often
derived from a class of domains $\cA$ in $\bbR^d$, in the following way:
%
%
\begin{equation} \label{eq:cluster-class}
\Km= \{K = A \cap\V_m\dvtx A \in\cA\}.
\end{equation}
Most of the literature assumes that the class $\cA$ is parametric,
exemplified by deformable templates, for which theoretical results are
available, especially in the case of the square lattice
\cite{boxscan,MGD,morel,jiang,perone,boutsikas}.
In particular, with a normal location model, the scan statistic
performs well, in the sense that it is asymptotically minimax; this is
shown in \cite{MGD} in a slightly different context tailored to image
processing applications. We also mention the recent work \cite{HJ09},
which considers the detection of multiple clusters (intervals) of
various amplitudes in the one-dimensional lattice.
As for nonparametric classes of domains, \cite{MGD} argues that the
scan statistic is asymptotically minimax for the case of star-shaped
clusters with smooth boundaries.

When $\V_m$ is endowed with a graph structure, \cite{maze} considers
paths of a certain length. In this setting, the scan statistic is shown
to be asymptotically minimax when the graph $\V_m$ is a complete,
regular tree and near-minimax for many other types of graphs, such as
the $d$-dimensional lattice for $d \geq3$. Addario-Berry et al. \cite{lugosi}
considers the same general testing problem with a focus on cluster
classes defined within the complete graph, such as cliques and spanning trees.
Note that part of the material presented here appeared in \cite{isicluster}.


\subsection{New theoretical results}

We describe here in an informal way the results we obtain.

In Section \ref{sec:geometric}, we focus on situations where the vertex
set $\V_m$ is embedded in a Euclidean space and well spread out in a
compact domain. Within this framework, we consider in Section
\ref{sec:thick} a geometric class of clusters obtained as in (\ref
{eq:cluster-class}) with $\cA$ a class of blobs that are mild
deformations of the unit ball. The clusters obtained in this way are
``thick,'' in the sense that they are not filamentary.
See Figure \ref{fig:thick}.
In particular, this class contains all the common parametric classes
obtained from parametric shapes such as hyperrectangles and ellipsoids,
as long as the shape is not too narrow. Note that the size, the (exact)
shape and the spatial location of the anomalous cluster under the
alternative is unknown.
In Corollary \ref{cor:allK}, we show that (under specific conditions)
$\hyp_0^m$ and $\hyp_1^m$ are asymptotically inseparable if there is
$\eta_m \to0$ slowly enough such that, for all $K \in\Km$,
\[
\Lambda_K \leq(1 - \eta_m) \sqrt{2 \log(m/|K|)};
\]
and conversely, we show that a version of the scan statistic
asymptotically separates $\hyp_0^m$ and $\hyp_1^m$ if there is $\eta_m
\to0$ slowly enough such that, for all $K \in\Km$,
\[
\Lambda_K \geq(1 + \eta_m) \sqrt{2 \log(m/|K|)}.
\]
Note that the detection rate is the same as for the class of balls so
that, perhaps surprisingly, scanning for the location (and not the
shape) is what drives the minimax detection risk.

In Section \ref{sec:thin}, we consider ``thin'' clusters, obtained as
in (\ref{eq:cluster-class}) with $\cA$ a class of ``bands'' around
smooth curves, surfaces or higher-dimensional submanifolds. In
particular, this class contains hyperrectangles and ellipsoids that are
sufficiently thin; see Figure \ref{fig:thin}.
It turns out that, contrary to what happens for thick clusters,
scanning for the actual shape impacts the minimax detection risk and
is, in fact, the main contributor for some nonparametric classes. The
situation is mathematically more challenging, yet we are able to prove
the following in Proposition \ref{prp:hold-lb}. Consider the class of
bands of thickness $r_m$ around $C^2$ curves of bounded curvature. Then
(under specific conditions), $\hyp_0^m$ and $\hyp_1^m$ are
asymptotically inseparable if
\[
\Lambdamax_m r_m^{1/4} (\log m)^{3/2} \to0.
\]
In Theorem \ref{thm:thin}, we show that, in the same setting, some
$\eps
_m$-scan statistic asymptotically separates $\hyp_0^m$ and $\hyp_1^m$ if
\[
\Lambdamin_m r_m^{1/4} \to\infty.
\]
Hence, some form of scan statistics achieves a detection rate within a
factor of $(\log m)^{3/2}$ from the minimax rate.

In Section \ref{sec:time}, we consider the spatio-temporal setting.
We first consider cluster sequences that admit a ``thick'' limit.
Cellular automata, which have been used to model epidemics
\cite{MR1683061}, satisfy this condition in some cases. In Proposition
\ref{prp:thick-time-ub}, we show that scanning over space--time cylinders,
as done in disease outbreak detection, achieves the asymptotic minimax risk.
We then consider cluster sequences with controlled space--time
variations, which may be a relevant model for applications such as
target tracking \cite{targetsensor}. We consider a fairly general
model in Proposition~\ref{prp:thin-time-ub}.

In Section \ref{sec:graph}, we assume that $\V_m = \{0, 1, \ldots,
m^{1/d}-1\}^d$, with $m^{1/d}$ an integer, seen as a subgraph of the
$d$-dimensional lattice. We first consider, in Section \ref{sec:band},
bands around nearest-neighbor paths;
see Figure \ref{fig:band}.
We extend the results obtained in \cite{maze} to paths. For example,
consider bands of thickness $h_m$ around a path of length $\ell_m$,
both powers of $m$. The bounds in Theorem \ref{thm:band} imply that
$\hyp_0^m$ and $\hyp_1^m$ are asymptotically inseparable if
\[
\Lambdamax_m (\ell_m/h_m)^{-1/2} (\log m)^{3/2} \to0.
\]
Conversely, Proposition \ref{prp:band-ub} states that an $\eps$-scan
statistic asymptotically separates $\hyp_0^m$ and $\hyp_1^m$ if
\[
\Lambdamin_m (\ell_m/h_m)^{-1/2} \to\infty.
\]
Therefore, some form of scan statistic is again within a factor of
$(\log m)^{3/2}$ from optimal.
In Section \ref{sec:animal}, we consider arbitrary connected
components, constraining only the size; see Figure \ref{fig:animal}.
In Proposition \ref{prp:animal}, we obtain a sharp detection rate for
clusters of very small size.

\subsection{Structure of the paper}

We have just described the contents of Sections~\ref{sec:geometric}
and \ref{sec:graph}. Section \ref{sec:discussion} is our discussion
section. We extend the results obtained for the normal location model
to any exponential family in Section \ref{sec:expo}. Other extensions
are described in Section \ref{sec:extensions}. We state some open
problems in Section \ref{sec:open}. In Section \ref{sec:computations},
we briefly discuss the challenge of computing the scan statistic. The
technical arguments are gathered in the \ref{supp}.

\subsection{Notation}

For two sequences of real numbers $(a_m)$ and $(b_m)$, $a_m \asymp b_m$
means that $a_m = O(b_m)$ and $b_m = O(a_m)$; $a_m \lesssim b_m$ means
that $a_m \leq(1 + o(1)) b_m$. For $a, b \in\bbR$, we use $a \vee b$
(resp., $a \wedge b$) to denote $\max(a,b)$ [resp., $\min(a,b)$]. For $a
\in\bbR$, let $[a]$ be the integer part of $a$; $\lfloor a \rfloor=
[a]$ if $a$ is not an integer and $[a]-1$ otherwise; and $\lceil a
\rceil= [a]+1$.
For a set $A$, $|A|$ denotes its cardinality.
Define $\log_\dag(x) = \log x$ if $x \geq e$ and $= 1$ otherwise.
All the limits in the text are when $m \to\infty$.
Throughout the paper, we use $C$ to denote a generic constant,
independent of $m$, whose particular value may change with each appearance.
We introduce additional notation in the text.

\section{Clusters as geometric shapes in Euclidean space}
\label{sec:geometric}

We assume that the nodes are embedded in $\Omega_d \subset\bbR^d$, a
compact set with nonempty interior. 
Let \mbox{$\| \cdot\|$} denote the corresponding Euclidean norm. For $A
\subset\Omega_d$ and $x \in\Omega_d$, let $\operatorname{dist}(x, A) =
{\inf_{y
\in A}} \|x - y\|$ and for $r > 0$, define
\[
B(A, r) = \{x \in\bbR^d\dvtx\operatorname{dist}(x, A) < r\}.
\]
In particular, $B(x, r)$ denotes the (open) Euclidean ball with center
$x$ and radius $r$. On occasion, we will add a subscript $d$ to
emphasize that this is a $d$-dimensional ball.

We consider a sequence $(\V_m)$ of finite subsets of $\Omega_d$, of
size $|\V_m| = m$, that are evenly spread out, in the following sense:
there is a constant $C \geq1$, independent of $m$ and a sequence
$r_m^* \to0$ such that
%
%
\begin{equation} \label{eq:unif}
C^{-1} m r^d \leq|B(x, r) \cap\V_m| \leq C m r^d\qquad \forall r
\in[r_m^*, 1], \forall x \in\Omega_d.
\end{equation}
In words, the number of nodes in any ball that is not too small is
roughly proportional to its volume.
For the regular lattice with $m$ nodes in $\Omega_d = [0,1]^d$,
condition (\ref{eq:unif}) is satisfied for $r_m^* > \sqrt{d} m^{-1/d}$.
This is the smallest possible order of magnitude; indeed, for some
constant $C > 0$ and $r$ small enough, there is a set with more than $C
r^{-d}$ disjoint balls with centers in $\Omega_d$, and, by (\ref
{eq:unif}), they are all nonempty if $r \geq r_m^*$, which forces
$r_m^* \geq C m^{-1/d}$.
Another example of interest is that of $\V_m$ obtained by sampling $m$
points from the uniform distribution, or any other distribution with a
density with respect to the Lebesgue measure on $\Omega_d$, bounded
away from zero and infinity; in that case, (\ref{eq:unif}) is satisfied
with high probability for $r_m^* \geq C (\log(m)/m)^{1/d}$ when $C$ is
large enough; for an extensive treatment of this situation, see
\cite{penrose}, Chapter 4.

\subsection{Thick clusters}
\label{sec:thick}

In this section, we consider clusters as in (\ref{eq:cluster-class}),
where $\cA$ is a class of bi-Lipschitz deformations of the unit
$d$-dimensional ball. This includes the vast majority of all the
parametric clusters considered in the literature, such as
hyperrectangles and ellipsoids, as long as the shape is not too narrow.
Note that a slightly less general situation is briefly mentioned
in \cite{MGD}.

We start with a lower bound on the minimax detection rate for discrete
balls of a given radius.
\begin{prp} \label{prp:balls}
Consider $\lambda_m \to0$ such that $\lambda_m \geq r_m^*$ and let
$\Km
$ be the class of all discrete balls of radius $\lambda_m$, that is,%
\[
\Km= \{K = B(x, \lambda_m) \cap\V_m \dvtx x \in\Omega_d\}.
\]
$\hyp_0^m$ and $\hyp_1^m$ are then asymptotically inseparable if
\[
\Lambdamax_m \leq\sqrt{2 d \log(1/\lambda_m)} -\eta_m,
\]
where $\eta_m \to\infty$.
\end{prp}


We now consider a much larger class of clusters and show that,
nevertheless, a form of scan statistic achieves that same detection rate.
For a function $f\dvtx A \subset\bbR^p \to\bbR^d$, its Lipschitz constant
is defined as
\[
\lambda_f = \sup_{x \neq y} \frac{\|f(x) - f(y)\|}{\|x - y\|}.
\]
%
For $\kappa\geq1$, let $\cF_{d,d}(\kappa)$ be the subclass of
bi-Lipschitz functions $f\dvtx\oB(0,1) \subset\bbR^d \to\Omega_d$ such
that $\lambda_f \lambda_{f^{-1}} \leq\kappa$ or, equivalently,
%
%
\begin{equation} \label{eq:F}
\sup_{x \neq y} \frac{\|f(x) - f(y)\|}{\|x - y\|} \leq\kappa\inf_{x
\neq y} \frac{\|f(x) - f(y)\|}{\|x - y\|}.
\end{equation}
%

For a function $f\dvtx A \to\bbR^d$, define
\[
K_f = \im(f) \cap\V_m,\qquad \im(f) := \{f(x)\dvtx x \in A\}.
\]
Note that $\lambda_f$ is intimately related to the size of $\im(f)$ and
therefore of $K_f$. Indeed, a simple application of (\ref{eq:F})
implies that, for any $f \in\cF_{d,d}(\kappa)$,
%
%
\begin{equation} \label{eq:imf}
B\bigl(f(0), \lambda_f/\kappa\bigr) \subset\im(f) \subset B(f(0), \lambda_f).
\end{equation}
This implies that sets of the form $\im(f)$, with $f \in\cF
_{d,d}(\kappa)$, are ``thick,'' in the sense that the smallest ball(s)
containing $\im(f)$ and the largest ball(s) included in $\im(f)$ are of
comparable sizes.
\begin{theorem} \label{thm:thick}
Consider $\lambda_m \to0$ such that $\lambda_m \geq r_m^*$ and define
\[
\Km= \{K_f\dvtx f \in\cF_{d,d}(\kappa), \lambda_f \geq\lambda_m\}.
\]
An $\eps_m$-scan statistic with $\eps_m \to0$ and $\eps_m (\log
(1/\lambda_m))^{1/(2d)} \to\infty$ then asymptotically separates
$\hyp
_0^m$ and $\hyp_1^m$ if
\[
\Lambdamin_m \geq\sqrt{2 d \log(1/\lambda_m)} +\eta_m,
\]
where $\eta_m = \eps_m^2 \sqrt{2 d \log(1/\lambda_m)}$.
Moreover, if $r_m^* \asymp m^{-1/d}$ and
\[
\Km= \{K_f\dvtx f \in\cF_{d,d}(\kappa)\},
\]
then an $\eps_m$-scan statistic with $\eps_m \to0$ and $\eps_m
(\log
m)^{1/(2d)} \to\infty$ asymptotically separates $\hyp_0^m$ and $\hyp
_1^m$ if
\[
\Lambdamin_m \geq\sqrt{2 \log m} +\eta_m,
\]
where $\eta_m = \eps_m^2 \sqrt{2 \log m}$.
\end{theorem}

Therefore, on a larger class of mild deformations of the unit ball,
some form of scan statistic achieves essentially the same detection
rate as for the class of balls stated in Proposition \ref{prp:balls}.

We note that the lower bound on $\Lambdamin_m$ is driven by the smaller
clusters in the class and that the performance guarantee is subject to
a proper choice of $\eps_m$. A~simple fix for both issues is to combine
the tests for different cluster sizes with an appropriate correction
for multiple testing. We summarize the consequence of Proposition
\ref{prp:balls} and Theorem \ref{thm:thick} with this observation in the
following result, inspired by \cite{boxscan}.
%
\begin{cor} \label{cor:allK}
Consider $\lambda_m \to0$ such that $\lambda_m \geq r_m^*$ and define
\[
\Km= \{K_f\dvtx f \in\cF_{d,d}(\kappa), \lambda_m \geq\lambda_f \geq
r_m^*\}.
\]
$\hyp_0^m$ and $\hyp_1^m$ are then asymptotically inseparable if, for
all $K \in\Km$,
\[
\Lambda_K \leq\sqrt{2 \log(m/|K|)} -\eta_m,
\]
where $\eta_m \to\infty$.
Conversely, let $T_\ell$ be an $\eps_{\ell}$-scan statistic for the
subclass $\{K_f \in\Km\dvtx2^{-\ell} \leq\lambda_f < 2^{-\ell+1}\}$ with
$\eps_{\ell} \ell^{1/(2d)} \to\infty$. There is a test based on $\{
T_\ell\dvtx\ell\geq0\}$ that asymptotically separates $\hyp_0^m$ and
$\hyp_1^m$ if, for all $K \in\Km$,
\[
\Lambda_K \geq\sqrt{2 \log(m/|K|)} +\eta_K,
\]
where $\eta_K = \eps_{\ell_K}^2 \sqrt{2 \log(m/|K|)}$ and $\ell_K =
\log(m/|K|)$.
\end{cor}

The same procedure, that is, combining $\eps$-scan statistics at
different (dyadic) scales, may be implemented in any of the settings we
consider in this paper to obtain a test that does not depend on a
tuning parameter like $\eps_m$ and achieves the same optimal rate at
every size. This is simply due to the fact that we only need to
consider the order of $\log m$ scales and the fast decaying tails of
the scan statistics under the null.


\subsubsection*{Union of thick clusters}
In a number of situations, the signal to be detected may be composed of
several clusters. Our results extend readily to this case.
Let $j_m$ be a positive integer and consider sets of the form $\bigcup
_{j = 1}^{j_m} K_{f_j}$, where the union is over some $f_j \in\cF
_{d,d}(\kappa)$ such that, for $j, j'$, $\lambda_{f_j} \leq C \lambda
_{f_{j'}}$ and
\[
\|f_j(0) - f_{j'}(0)\| \leq C (\lambda_{f_{j}} \vee\lambda_{f_{j'}}),
\]
so that the sets $\im(f_j)$ and $\im(f_{j'})$ are of comparable sizes
and not too far from each other. In that case, Theorem \ref{thm:thick}
applies unchanged, as long as the number of clusters is not too large,
specifically if $j_m = o(\log(1/\lambda_m))^{1/d}$.
(This can be improved if the $K_{f_j}$'s do not overlap too much.)
If the proximity constraint is dropped, then the term $\log(1/\lambda
_m)$ in Theorem \ref{thm:thick} is replaced by $j_m \log(1/\lambda_m)$.

\subsection{Thin clusters}
\label{sec:thin}

In this section, we consider clusters that are built from smooth
embeddings in $\Omega_d$ of the unit $p$-dimensional ball, where $p <
d$. The special case of curves ($p = 1$) is, for example, relevant in
road tracking \cite{roadtracking} and in modeling blood vessels in
medical imaging \cite{bloodvessel04}. As in the previous section, the
results we obtain below are valid for (some) unions of such subsets
and, in particular, for submanifolds with a wide array of topologies.

For a differentiable function $f$ between two Euclidean spaces, let
$Df$ denote its Jacobian matrix. For $\kappa\geq1$, let $\cF
_{p,d}(\kappa)$ be the class of twice differentiable, one-to-one
functions $f\dvtx\oB(0,1) \subset\bbR^p \to\im(f) \subset\Omega_d$
satisfying $\lambda_f \lambda_{f^{-1}} \leq\kappa$ and $\lambda_{Df}
\leq\kappa\lambda_f$.
We consider clusters that are tubular regions around the range of
functions in $\cF_{p,d}(\kappa)$.
For a function $f$ with values in $\bbR^d$ and $r > 0$, define
\[
K_{f,r} = B(\im(f), r) \cap\V_m.
\]
%
Again, $\lambda_f$ is intimately related to the size of $B(\im(f),r)$
and $K_{f,r}$. This relationship is made explicit in the \ref{supp}.
We consider classes of clusters of the form $\{K_{f,r}\dvtx f \in\cF\}$,
where $\cF$ is a subclass of $\cF_{p,d}(\kappa)$.

We start with a result on the performance of the scan statistic.
For a class $\cF$ of functions with values in $\bbR^d$ and for $\eps>
0$, let $N_\eps(\cF)$ denote its $\eps$-covering number for the
sup-norm, that is,
\[
N_\eps(\cF) = \min\Bigl\{n\dvtx\exists f_1, \ldots, f_n \in\cF\mbox{,
s.t. }
{\max_{f \in\cF} \min_j} \|f - f_j\|_\infty\leq\eps\Bigr\}.
\]
%
%
\begin{theorem}
\label{thm:thin}
Let $C$ be the constant defined in Lemma \textup{B.2} in the
\setattribute{ref}{fmt}\textit\ref{supp}\setattribute{ref}{fmt}\textup.
Consider $\lambda_m, r_m \to0$ such that $C^{-1} \lambda_m \geq r_m
\geq r_m^*$ and let $\cF$ be a subclass of $\cF_{p,d}(\kappa)$.
Define
\[
\Km= \{K_{f,r}\dvtx f \in\cF, \lambda_f \geq\lambda_m, C^{-1} \lambda_m
\geq r \geq r_m\}.
\]
An $\eps_m$-scan statistic with $\eps_m = o(r_m^{1/2})$ then
asymptotically separates $\hyp_0^m$ and $\hyp_1^m$ if
\[
\Lambdamin_m \geq(1 +\eps_m^2) \sqrt{2 \log N_{\eps_m^2}(\cF) + 2 d
\log(1/\lambda_m)}.
\]
\end{theorem}

Just as in Theorem \ref{thm:thick}, if $r_m^* \asymp m^{-1/d}$, we can
dispense with the restriction $r_m \geq r_m^*$ and replace the factor
$\log(1/\lambda_m^d)$ by $\log m$ in the bound.

For a typical parametric class $\cF$, $\log N_\eps(\cF) \sim a(\cF)
\log(1/\eps)$, so the scan statistic (over an appropriate net) is
accurate if
%
%
\begin{equation} \label{eq:thin-p}
\Lambdamin_m \geq(1 +r_m) \sqrt{2 a(\cF) \log(1/r_m) + 2 d \log
(1/\lambda_m)}.
\end{equation}
On the other hand, $\log N_\eps(\cF) \asymp(1/\eps)^{a(\cF)}$ for a
typical nonparametric class $\cF$ \cite{MR0124720}, so the scan
statistic (over an appropriate net) is accurate if
%
%
\begin{equation} \label{eq:thin-np}
\Lambdamin_m r_m^{a(\cF)/2} \to\infty.
\end{equation}

Finding a sharp lower bound for the minimax detection rate is more
challenging for thin clusters compared to thick clusters. By
considering disjoint tubes around $p$-dimensional hyperrectangles, we
obtain a lower bound that matches, in order of magnitude, the rate
achieved by the scan statistic when the class $\cF$ is parametric,
displayed in (\ref{eq:thin-p}).
%
\begin{prp} \label{prp:thin-lb}
Consider $\lambda_m, r_m \to0$ with $\lambda_m \geq r_m \geq r_m^*$.
Let $U\dvtx\bbR^p \to\bbR^d$ be the canonical embedding so that $U x =
(x,0)$ and let
\[
\cF= \{f\dvtx\oB(0,1) \subset\bbR^p \to\Omega_d, f(x) = \lambda_m
U x
+ b\mbox{,  where } b \in\bbR^d\}.
\]
Define
\[
\Km= \{K_{f,r_m}\dvtx f \in\cF\}.
\]
$\hyp_0^m$ and $\hyp_1^m$ are then asymptotically inseparable if
\[
\Lambdamin_m \leq\sqrt{2 (d-p) \log(1/r_m) + 2 p \log(1/\lambda_m)}
-\eta_m,
\]
where $\eta_m \to\infty$.
\end{prp}

The proof is parallel to that of Proposition \ref{prp:balls} and is
therefore omitted.

For at least one family of nonparametric curves ($p = 1$), we show that
the rate displayed at (\ref{eq:thin-np}) matches the minimax rate,
except for a logarithmic factor.
For concreteness, we assume that $\Omega_d = [0,1]^d$.
Let $\cH(\alpha, \kappa)$ be the H\"older class of functions $g\dvtx[0,1]
\to[0,1]$ satisfying
\begin{eqnarray*}
\bigl|g^{(s)}(x)\bigr| &\leq& \kappa\qquad \forall x \in[0,1], \forall s < \alpha;
\\
\bigl|g^{(\lfloor\alpha\rfloor)}(x) - g^{(\lfloor\alpha\rfloor)}(y)\bigr|
&\leq& \kappa|x - y|^{\alpha- \lfloor\alpha\rfloor}\qquad \forall x, y
\in[0,1].
\end{eqnarray*}
\begin{prp} \label{prp:hold-lb}
Let $r_m \to0$ with $r_m \geq r_m^*$.
Let $\cF$ be the class of functions of the form $f(x) = (x, g_1(x),
\ldots, g_{d-1}(x))$, where $g_j \in\cH(\alpha, \kappa)$, with
$\alpha
\geq2$.
Define
\[
\Km= \{K_{f,r_m}\dvtx f \in\cF\}.
\]
$\hyp_0^m$ and $\hyp_1^m$ are then asymptotically inseparable if
\[
\Lambdamax_m r_m^{1/(2\alpha)} (\log m)^{3/2} \to0.
\]
%
\end{prp}

Thus, for the detection of curves with H\"older regularity, a scan
statistic achieves the minimax rate within a poly-logarithmic factor.
We prove Proposition \ref{prp:hold-lb} by reducing the problem of
detecting a band in a graph so that we can use results from Section
\ref{sec:band}. We do not know how to generalize this approach to
higher-dimensional surfaces (i.e., $p \geq2$).

\subsection{The spatio-temporal setting}
\label{sec:time}

In this section, we consider the spatio-temporal setting described in
Section \ref{sec:dynamic}. This is a special case of the spatial
setting we have considered thus far, with time playing the role of an
additional dimension.
For their relevance in applications and concreteness of exposition, we
focus on two specific models. In Section \ref{sec:thick-time}, we
consider cluster sequences with a limit; as we shall see, this
assumption is implicit in some popular models for epidemics. In Section
\ref{sec:thin-time}, we consider cluster sequences of bounded variations.

In the remainder of this section, we assume, for concreteness, that $\T
_m = \{0, 1, \ldots, t_m\}$ with $t_m \to\infty$.
Our results apply without any changes if the set of nodes varies with
time, that is, with index set of the form ${\prod_{t \in\T_m}} \V
_m^{t}$, in the case where each $\V_m^{t}$ satisfies (\ref{eq:unif})
with $C$ and $r_m^*$ independent of $t$.

\subsubsection{Cluster sequences with a limit}
\label{sec:thick-time}

We focus here on cluster sequences obeying $K_{t_m} \neq\varnothing$,
that is, the anomalous cluster is present at the last time point. This
is a standing assumption in syndromic surveillance systems
\cite{diseaseoutbreak}. To illustrate the difference, consider a typical
change-point problem setting, where $\V_m$ contains only one node and,
for simplicity, assume that $\Lambda_K$ is independent of $K$ and that
$\Lambda_m$ denotes this common value. First, let the cluster be any
discrete interval (in time), so the signal may not be present at time
$t = t_m$. This is a special case of Section \ref{sec:thick}, with time
playing the role of a spatial dimension ($d = 1$); we saw in Corollary
\ref{cor:allK} that the detection threshold is at $\Lambda_m \sim
\sqrt
{2 \log|\T_m|}$. Now, let the emerging cluster be any discrete
interval that includes $t = t_m$. Detecting such an interval is
actually much easier since we do not need to search where the interval
is located, which is what drives the detection threshold for the thick
clusters in Section \ref{sec:thick}---we need only determine its
length. Specifically, the scan statistic over the dyadic intervals
containing $t = t_m$ asymptotically separates the hypotheses if
$\Lambda
_m \asymp\sqrt{{\log\log}|\T_m|}$.

Regarding the actual evolution of the cluster in time, a number of
growth models have been suggested, for example, cellular automata
\cite{MR2367301,MR1849342} and their random equivalent, threshold growth
automata \cite{MR2206346,MR1698409}, which have been used to model
epidemics \cite{MR1683061}. The latter includes the well-known
Richardson model \cite{Richardson1973jk}. Under some conditions, these
models develop an asymptotic shape (with probability one), a convex
polygon in the case of threshold growth automata. Less relevant for
modeling epidemics, internal diffusion limited aggregation is another
growth model with a limiting shape \cite{Lawler1992jk}.

The simplest cluster sequences with limiting shape are space--time
cylinders, for which we have the equivalent of Proposition
\ref{prp:balls}. (The proof is completely parallel and we omit details.)
\begin{prp} \label{prp:thick-time-lb}
Consider $\lambda_m \to0$ with $\lambda_m \geq r_m^*$ and let $\Km$ be
the class of all space--time cylinders of the form $K_t = B(x, \lambda
_m) \cap\V_m, \forall t = 0, \ldots, t_m$, where $x \in\Omega_d$.
$\hyp_0^m$ and $\hyp_1^m$ are then asymptotically inseparable if
\[
\Lambdamin_m \leq\sqrt{2 d \log(1/\lambda_m)} -\eta_m,
\]
where $\eta_m \to\infty$.
\end{prp}

With only one possible shape and known starting point, such a model is
rather uninteresting. We now consider a much larger class of cluster
sequences with some sort of limit [in the sense of (\ref
{eq:limit-set})] and show that, nevertheless, a form of scan statistic
achieves that same detection rate. For a cluster sequence $K = (K_t, t
\in\T_m)$, let $t_{\Kt} = \min\{t\dvtx K_t \neq\varnothing\}$, which is the
time when $K$ originates.
The following is the equivalent of Theorem \ref{thm:thick}.
The metric $\delta$ appearing below is defined in (\ref{eq:delta}).
\begin{prp} \label{prp:thick-time-ub}
Consider sequences $\lambda_m \to0$ with $\lambda_m \geq r_m^*$ and
$\log\log t_m = o(\log(1/\lambda_m))$,
and a function $\nu(t)$ with $\lim_{t \to\infty} \nu(t) = 0$ and
$\nu
(t) \leq1$ for all $t \geq0$. Let $\Km$ be a class of cluster
sequences such that $t_m - \max\{t_K \dvtx K \in\Km\} \to\infty$ and, for
each $K = (K_t, t \in\T_m) \in\Km$, there exists $f \in\cF
_{d,d}(\kappa)$ with $\lambda_f \geq\lambda_m$ such that
%
%
\begin{equation} \label{eq:limit-set}
\delta\bigl(K_{t}, \im(f) \cap\V_m\bigr) \leq\nu(t-t_{\Kt}) \qquad
\forall t \in \T_m.
\end{equation}
There is then a scan statistic over a family of space--time cylinders
that asymptotically separates $\hyp_0^m$ and $\hyp_1^m$ if
\[
\Lambdamin_m \geq(1 + \xi_m) \sqrt{2 d \log(1/\lambda_m)},
\]
where $\xi_m \to0$ slowly enough.
\end{prp}

If the starting time is uniformly bounded away from $t_m$ and the
convergence to the thick spatial cluster [in the sense of (\ref
{eq:limit-set})] occurs at a uniform speed, then all of the cluster
sequences in the class have sufficient time to develop into their
``limiting'' shapes. The space--time cylinders over which we scan are
based on an $\eps$-net for the possible limiting shapes, that is, the
class of thick clusters.

Scanning over space--time cylinders (with balls as bases) is advocated
in the disease outbreak detection literature \cite{diseaseoutbreak}.
Although seemingly naive, this approach achieves, in our asymptotic
setting, the minimax detection rate if the cluster sequences develop
into balls and, in general, falls short by a constant factor.

We mention that the equivalent of Corollary \ref{cor:allK} holds here
as well, in that we can combine the different scans at different
space--time scales to obtain a test that does not depend on a tuning
parameter (implicit here) and which achieves the same rate for the
cluster class defined as above, but with $\lambda_m \geq\lambda_f
\geq
r_m^*$, which is the class that appears in Corollary \ref{cor:allK}.

\subsubsection{Cluster sequences of bounded variation}
\label{sec:thin-time}

In target tracking \cite{linesand,targetsensor}, the target is
usually assumed to be limited in its movements due to maximum speed and
maneuverability. With this example in mind, we consider classes of
cluster sequences of bounded variation, meaning that the cluster is
limited in the amount it can change in a given period of time. As the
rates we obtain in this subsection are the same with or without the
condition $K_{t_m} \neq\varnothing$, we do not make that assumption. Let
$t_K^+ = \max\{t\dvtx K_t \neq\varnothing\}$.

We consider space--time tubes around H\"older space--time curves. For
$\alpha\in(0,1]$ and $\kappa> 0$, let $\cH_\infty(\alpha, \kappa)$
be the H\"older class of functions $g\dvtx[0,\infty) \to[0,1]$ satisfying
%
%
\begin{equation} \label{eq:Halpha}
|g(x) - g(y)| \leq\kappa|x - y|^{\alpha} \qquad\forall x, y \in[0,
\infty).
\end{equation}
The following is the equivalent of Proposition \ref{prp:hold-lb}.
\begin{prp} \label{prp:thin-time}
Assume that $\Omega_d = [0,1]^d$.
Consider sequences $r_m \to0$ with $r_m \geq2 r_m^*$ and $\xi_m$ such
that $1 \leq\xi_m \leq t_m$. Let $\Km$ be the class of all cluster
sequences $K_g$ of the form $K_{g,t} = B(g(t/\xi_m), r_m) \cap\V_m$
for all $t = t_K, \ldots, t_K^+$, for some $g = (g_1, \ldots, g_d)$ with
$g_j \in\cH_\infty(\alpha, \kappa)$.
Then, $\hyp_0^m$ and $\hyp_1^m$ are asymptotically inseparable if
\[
\Lambdamax_m (t_m/\lceil\xi_m r_m^{1/\alpha} \rceil)^{-1/2} \log
(t_m/\lceil\xi_m r_m^{1/\alpha} \rceil) \bigl(\log(\xi_m) + \log\log
(t_m)\bigr)^{1/2} \to0.
\]
Conversely, an $\eps$-scan statistic with $\eps< \sqrt{2}$
asymptotically separates $\hyp_0^m$ and $\hyp_1^m$ if
\[
\Lambdamin_m \bigl((t_m/\lceil\xi_m r_m^{1/\alpha} \rceil) \log_\dag
(\xi
_m^{-1} r_m^{-1/\alpha}) + \log m\bigr)^{-1/2} \to\infty.
\]
\end{prp}

For simplicity, assume that $t_m$ is a power of $m$. If $\xi_m
r_m^{1/\alpha} = O(1)$, then the detection threshold is roughly of
order $t_m^{1/2}$, while if $\xi_m r_m^{1/\alpha}$ is large, yet small
enough that $t_m/(\xi_m r_m^{1/\alpha})$ is still a power of $m$, then
the detection threshold is roughly of order $(t_m/(\xi_m r_m^{1/\alpha
}))^{1/2}$.

A form of scan statistic is actually able to attain the same detection
rate when the radius is unknown, but restricted to $r \geq r_m$. In
fact, another form of scan statistic achieves a slightly different rate
over a much larger class of cluster sequences with bounded variations.
Let $\cS(r, \kappa)$ be the set of subsets $S \subset\Omega_d$ such
that $B(x, r) \subset S \subset B(x, \kappa r)$ for some $x \in\Omega_d$.
\begin{prp} \label{prp:thin-time-ub}
Consider a sequence $\xi_m$ such that $1 \leq\xi_m \leq t_m$ and a
constant $\eta> 0$. Define $\Km$ as the class of cluster sequences $K$
such that, for each $t = t_K, \ldots, t_K^+$, $K_t = S_t \cap\V_m$,
where $S_t \in\cS(r_t, \kappa)$ for some $r_t \geq r_m^*$, and, for
any $s, t = t_K, \ldots, t_K^+$,
%
%
\begin{equation} \label{eq:delta-eta-xi}
\delta(K_t, K_s) \leq\eta \qquad\mbox{if } |t - s| \leq\xi_m.
\end{equation}
%
Then, for $\eta$ small enough, an $\eps$-scan statistic with $\eps<
\sqrt{2}$ asymptotically separates $\hyp_0^m$ and $\hyp_1^m$ if
\[
\Lambdamin_m \bigl((t_m/\xi_m) (\log m) + \log t_m\bigr)^{-1/2} \to\infty.
\]
\end{prp}

Consider the condition
%
%
\begin{equation} \label{eq:delta-nu}
\delta(K_t, K_s) \leq\nu(|t-s|/\xi_m) \qquad\forall s, t \in\{t_K,
\ldots, t_K^+\},
\end{equation}
for a function $\nu\dvtx[0, \infty) \to[0, \sqrt{2}]$. Then, (\ref
{eq:delta-eta-xi}) is satisfied with $\eta= \nu(1)$ and the same $\xi
_m$. The requirement in Proposition \ref{prp:thin-time-ub} is that
$\nu
(1)$ be small enough.
In particular, the cluster sequences considered in Proposition
\ref{prp:thin-time} satisfy, for some constant $C > 0$,
\[
\delta(K_t, K_s) \leq C r_m^{-1/2} \bigl(r_m^* \vee(|t-s|/\xi
_m)^\alpha\bigr)^{1/2} \qquad\forall s, t \in\{t_K, \ldots, t_K^+\}.
\]
This comes from Lemma C.1 in the \ref{supp} and (\ref{eq:Halpha}).
Therefore, assuming $\xi_m \ll(r_m^*)^{-1/\alpha}$, (\ref
{eq:delta-nu}) is satisfied with $\nu(u) = u^{\alpha/2}$ and $\xi_m$
replaced by $\xi_m r_m^{1/\alpha}$. In that case, the detection rates
obtained by the scan statistics of Propositions \ref{prp:thin-time} and
\ref{prp:thin-time-ub} are of comparable orders of magnitude.

\section{Clusters as connected components in a graph}
\label{sec:graph}

In this section, we model the network with the $d$-dimensional square
lattice; specifically, we assume that $m^{1/d}$ is an integer (for
convenience) and consider $\V_m = \{0, 1, \ldots,\break m^{1/d}-1\}^d$, seen
as a subgraph of the usual $d$-dimensional lattice. We assume that $d
\geq2$ since the case where $d = 1$ is treated in Section \ref{sec:thick}.
We work with the $\ell^1$-norm, which corresponds to the shortest-path
distance in the graph; let $B(v, h)$ denote the corresponding open ball
with center $v$ and radius $h$ so that $B(v,h) = \{v\}$ for $h \in
(0,1]$, and, for a subset of nodes $V$, let $B(V, h) = \bigcup_{v \in
V} B(v, h)$.

\subsection{Paths and bands}
\label{sec:band}

A nearest-neighbor band of length $\ell$ and width $h$ is of the form
$B(V, h)$, where $V = (v_0, \ldots, v_{\ell})$ forms a path in $\Zd$. A
band with unit width ($h = 1$) is just a path.

We say that a path $(v_0, \ldots, v_{\ell})$ in $\bbZ^d$ is
nondecreasing if, for all $t = 1, \ldots, \ell$, $v_t - v_{t-1}$ has
exactly one coordinate equal to 1 and all other coordinates equal to~0.
The case of paths was treated in detail in \cite{maze}; it corresponds
to taking $h_m = 1$ below.
\begin{theorem} \label{thm:band}
Suppose that $d \geq2$ and let $\Km$ be the class of bands of width
$h_m$ generated by nondecreasing paths in $\V_m$ of length $\ell_m$,
starting at the origin, with $m^{1/d} \geq\ell_m \geq h_m$.
Then, $\hyp_0^m$ and $\hyp_1^m$ are asymptotically inseparable if
\begin{eqnarray*}
\Lambdamax_m (\ell_m/h_m)^{-1/2} \log_\dag(\ell_m) (\log h_m +
\log_\dag
\log\ell_m)^{1/2} &\to& 0\qquad \mbox{for } d = 2, \\
\Lambdamax_m (\ell_m/h_m)^{-1/2} (\log_\dag h_m) (\log_\dag\log h_m)
&\to& 0\qquad \mbox{for } d \geq3.
\end{eqnarray*}
Conversely, an $\eps$-scan statistic with $\eps< \sqrt{2}$ fixed
asymptotically separates $\hyp_0^m$ and $\hyp_1^m$ if
\[
\Lambdamin_m (\ell_m/h_m)^{-1/2} \to\infty.
\]
\end{theorem}

For the case of nondecreasing paths, a form of the scan statistic
achieves the minimax rate in dimension $d \geq3$, while it falls short
by a logarithmic factor in dimension $d = 2$. In the latter
setting, Arias-Castro et al. \cite{maze} introduces a test that asymptotically separates
$\hyp_0^m$ and $\hyp_1^m$ if
%
%
\begin{equation} \label{eq:2d-lattice-ub}
\Lambdamin_m (\ell_m/h_m)^{-1/2} \log_\dag(\ell_m)^{1/2} \to
\infty,
\end{equation}
coming slightly closer to the minimax rate.

In fact, even when the band has unknown length, width and starting
location, and when the path is not restricted to be nondecreasing, a
form of scan statistic achieves the same rate, except for a logarithmic factor.
\begin{prp} \label{prp:band-ub}
Suppose that $d \geq2$ and let $\cK_m$ be the class of all bands of
width $h$ and length $\ell$, where $\ell_m \geq\ell\geq h \geq h_m$,
that are within $\V_m$ and generated by paths that do not self-intersect.
An $\eps$-scan statistic, with $\eps< \sqrt{2}$, then asymptotically
separates $\hyp_0^m$ and $\hyp_1^m$ if
\[
\Lambdamin_m \bigl(\ell_m/h_m + \log(m/h_m^d) + \log_\dag\log\ell
_m\bigr)^{-1/2} \to\infty.
\]
\end{prp}


\subsection{Arbitrary connected components}
\label{sec:animal}

We consider here classes of connected components with a constraint on
their sizes.
Arbitrary connected components in the square lattice are sometimes
called \textit{animals} or \textit{polyominoes} (polycubes in dimension $d
\geq3$), which are well-studied objects in combinatorics, where the
goal is to count the number of polyominoes \cite{Klarner1967fp}.
We mention in passing the results in \cite{MR1825148} which provide a
law of large numbers for the scan statistic under the null.
Otherwise, such objects are fairly new to statistics.
Detecting animals is, of course, harder than detecting paths since
paths are themselves animals. The result below offers a sharp detection
threshold for connected components of sufficiently small size.
\begin{prp} \label{prp:animal}
Let $\Km$ be the class of animals of size $k_m = o(m)$ within~$\V_m$.
$\hyp_0^m$ and $\hyp_1^m$ are then asymptotically inseparable if
\[
\Lambdamax_m \leq\sqrt{2 \log m} -\eta_m,
\]
where $\eta_m \to\infty$.
Conversely, let $\Km$ be the class of animals of size not exceeding
$k^+_m = o(\log m)$ within $\V_m$. The actual scan statistic then
asymptotically separates $\hyp_0^m$ and $\hyp_1^m$ if
\[
\Lambdamin_m \geq\sqrt{2 \log m}.
\]
\end{prp}

Note that, in general, we can obtain a quick (naive) upper bound on the
detection rate for large clusters by considering the simple test that
rejects for large values of $\sum_{v \in\V_m} X_v$ (this is the
``average test'' in \cite{lugosi}). This test asymptotically separates
$\hyp_0^m$ and $\hyp_1^m$ if $\Lambdamin_m k_m^{1/2} m^{-1/2} \to
\infty
$, assuming the clusters in $\Km$ are of size bounded below by $k_m$.
An open question of theoretical interest is whether, for the class of
animals of size $k_m = \sqrt{m}$ in the two-dimensional lattice, there
is a test that asymptotically separates $\hyp_0^m$ and $\hyp_1^m$ when
$\Lambdamin_m k_m^{-1/2} \to0$ slowly enough. In dimension three or
higher, Theorem \ref{thm:band} implies that this is not possible, even
for paths.

\section{Discussion}
\label{sec:discussion}

\subsection{Extension to exponential families}
\label{sec:expo}

Although the previous results were stated for the normal location
model, they extend to any one-parameter exponential model if the
anomalous clusters are large enough. For example, consider a Bernoulli
model where the variables are Bernoulli with parameter $1/2$ under the
null and with parameter $p_{K} > 1/2$ when they belong to the anomalous
cluster~$K$; or, a Poisson model where the variables are Poisson with
mean 1 under the null and $\mu_K > 1$ when they belong to the anomalous
cluster $K$. In general, transforming the variables and/or the
parameter if necessary, we may assume that the model is of the form
$F_\theta$, with density $f_\theta(x) = \exp(\theta x - \log
\ph(\theta))$ with respect to $F_0$, where, by definition,
$\ph(\theta) = \E_0 [\exp(\theta X)]$, where $\E_0$ denotes the
expectation under $F_0$.
We always assume that $\ph(\theta)<\infty$ for $\theta$ in a
neighborhood of $0$. Let $\sigma^2 = \var_0(X)$, the variance of $X
\sim F_0$. In the Bernoulli model, the correspondence is $\theta= \log
(p/(1-p))$ and $\sigma^2 = 1/4$; in the Poisson model, $\theta= \log
\lambda$ and $\sigma^2 = 1$.
Under the null hypothesis, all of the variables at the nodes have
distribution $F_0$, that is,
\[
\hyp_0^m \dvtx X_v \sim F_0 \qquad\forall v \in\V_m.
\]
Under the alternative, the variables at the nodes belonging to the
anomalous cluster $K \in\Km$ have distribution $F_{\theta_K}$ with
$\theta_{K} := \sigma\Lambda_K |K|^{-1/2}$, that is,
\[
\hyp_{1,K}^m \dvtx X_v \sim F_{\theta_K}\qquad \forall v \in K;\qquad X_v \sim
F_0\qquad
\forall v \notin K.
\]
As before, the variables are assumed to be independent.

If the clusters in the class are sufficient large, then the results
presented for the normal location family hold unchanged. Intuitively,
large enough clusters allow for the sums over them to be approximately
normally distributed. Details are provided in the \ref{supp}.
For example, we have the following equivalent of
Corollary~\ref{cor:allK} in the context of thick clusters as in Section
\ref{sec:thick}. Consider $\lambda_m \geq r_m \geq r_m^*$ with $m r_m^d
(\log1/r_m)^{-3} \to\infty$ (which guarantees that the clusters in
the class are large enough) and define the class
\[
\Km= \{K_f\dvtx f \in\cF_{d,d}(\kappa), \lambda_m \geq\lambda_f \geq
r_m\}.
\]
In this setting, under the Bernoulli model, the detection threshold is at
\[
p_{K} = \frac{1}{2} + \frac{1}{8 |K|^{1/2}} \bigl(2 \log(m/|K|)\bigr)^{1/2};
\]
under the Poisson model, the detection threshold is at
\[
\mu_K = 1 + \frac{1}{|K|^{1/2}} \bigl(2 \log(m/|K|)\bigr)^{1/2}.
\]

Note that without a lower bound on the minimum size of the anomalous
clusters, the general analysis breaks down and the results depend on
the specific exponential model. For example, unless $\min\{|K|\dvtx K
\in
\Km\} \to\infty$ fast enough, detection is impossible in the Bernoulli
model, even if the anomalous nodes have value 1 under the alternative.

\subsection{Other extensions}
\label{sec:extensions}
The array of possible models is as wide as the breadth of real-world
applications. We mention a few possible variations below.

\subsubsection*{Beyond exponential families}
Using an exponential family of distributions allows us to obtain sharp
detection lower bounds. Otherwise, similar results, although not as
sharp, may be obtained for essentially any family of
distribution~$F_\theta$, where the distance between the null $\theta= 0$ and an
alternative $\theta$ is in terms of the chi-square distance between
$F_0$ and $F_\theta$; see \cite{maze}, Section~5.

\subsubsection*{Different means at the nodes}
We could consider a situation where the mean varies over the nodes of
the anomalous cluster. This situation is considered in \cite{HJ09} for
the case of intervals, and the constant in the detection rate is indeed
different. We implicitly considered a worst case scenario where the
mean is bounded below over the anomalous cluster and subsequently
assumed it was equal to that lower bound everywhere over the anomalous
cluster. However, our results hold unchanged if we allow $X_v$ to have
any mean above $\theta_K$, for every $v \in K$, $K$ being the
anomalous cluster.

\subsubsection*{Dependencies}
Also of interest is the case where the variables are dependent. In the
spatial setting, the same paper \cite{HJ09} solves this problem for the
case of the one-dimensional lattice, with the correlation between $X_v$
and $X_w$ decaying as a function of distance between $v$ and $w$. We
postulate that the same result holds in higher dimensions. In the
spatio-temporal setting, variables could be dependent across time as
well, involving a higher degree of sophistication. We plan on pursuing
these generalizations in future publications.

\subsubsection*{Unknown variance or other parameters}
We assumed throughout that the variance was known (and equal to 1 after
normalization).\vadjust{\goodbreak} This is, in fact, a mild assumption, as one can
consistently estimate the variance using a robust estimator, say the
median absolute deviation (MAD), with the usual $\sqrt{m}$-convergence
rate, assuming that the anomalous cluster corresponds to a small part
of the entire network. When dealing with one-parameter families such as
Bernoulli or Poisson, the issue is to estimate the parameter under the
null and a robust version of the maximum likelihood (e.g., trimmed mean
for these two examples) can be used for that purpose.

\subsection{Energy, bandwidth and other constraints}

We assume throughout that a central processor has access to all the
information measured at the nodes and, based on that, makes a decision
as to whether there is an anomalous cluster of nodes in the network or
not. This assumption is reasonable in, for example, the context of
image processing or syndromic surveillance. However, real-world sensor
networks of the wireless type are often constrained by energy and/or
bandwidth considerations. A growing body of literature
\cite{targetsensor} is dedicated to designing efficient
(e.g., decentralized) communication protocols for sensor networks under
such constraints. As mentioned in Section \ref{sec:multi-hypo}, the
papers we are aware of consider very simplistic detection settings. In
the context of the present paper, it would be interesting to study how
the detection rates change when different communication protocols are used.

We also assume that we have infinite computational power. However, all
real-world systems operate under finite energy and processing
resources. In the same way, it would be interesting to know what
detection rates are achievable under such computational constraints.

\subsection{On computing the scan statistic}
\label{sec:computations}

In all of the settings we consider in this paper, the scan statistic
comes close to achieving the minimax detection rate. Turning to
computational issues, however, it is very demanding, even when scanning
for simple parametric clusters such as rectangles. For general shapes,
Duczmal, Kulldorff and Huang \cite{duczmal2006ess} suggests a simulated annealing
algorithm, which, from a theoretical point of view, is extremely
difficult to analyze.
For parametric shapes and blobs, Arias-Castro, Donoho and Huo \cite{MGD} advocates the use
of $\eps_m$-scan statistics based on multiscale nets built out of
unions of dyadic hypercubes; similar ideas appear in \cite{boxscan}.
Partial results suggest that this approach yields, in theory, a
near-optimal algorithm for detecting the more general thick clusters
considered in Section \ref{sec:thick}.

For the thin clusters of Section \ref{sec:thin}, or for the bands of
Section \ref{sec:band}, the situation is quite different. Take the
latter. After pre-processing the data by performing a moving average
with an appropriate radius, it remains to find the maximum over a
restricted, yet exponentially large, set of paths. Without further
restriction, this problem, known as the ``bank robber problem'' or
``reward budget problem'' \cite{DasGuptaHespanhaRiehlSontag06}, is
NP-hard. Note that DasGupta et al. \cite{DasGuptaHespanhaRiehlSontag06} suggests a
polynomial time approximation that deserves further investigation. The
case of thin clusters is even harder. In the context of point clouds,
Arias-Castro, Efros and Levi \cite{AriasCastro2009} introduces multiscale nets that could
be adapted to the setting of a network. It remains to compute the scan
statistic over this net, which seems particularly challenging for
surfaces of dimension $p \geq2$, which no longer correspond to paths.
In the spatio-temporal setting of Section \ref{sec:time}, dynamic
programming ideas could be used, as done in \cite{MSDFS} in the context
of point clouds and in \cite{chirpletpursuit} in the context of a
harmonic analysis decomposition of chirps.

\subsection{Open theoretical problems}
\label{sec:open}

The paper leaves two main theoretical problems unresolved.
The first one concerns obtaining sharper bounds for the detection of
thin clusters.
This is in the context of Section \ref{sec:thin}. For parametric
classes, the challenge is to match constants in the rate, while, for
nonparametric classes, the challenge is to obtain sharper lower bounds,
perhaps closer to what a scan statistic is shown to achieve in Theorem
\ref{thm:thick}. We were only able to do the latter for curves; see
Proposition \ref{prp:hold-lb}.

The second one concerns comparing the detection rates for arbitrary
connected components and for paths.
At a given size, the thicker the band (relative to its length), the
easier it is to detect it; see Theorem \ref{thm:band}. It seems,
therefore, that the most difficult connected components to detect are
paths or unions of paths. But is this true? In other words, are the
minimax detection rates for arbitrary connected components and paths of
a similar order of magnitude?

\section*{Acknowledgments}
The authors are grateful to the anonymous referees for suggesting an
expansion of the discussion section, for encouraging them to obtain
sharper bounds and for alerting them of the possibility of improving on
the performance of the scan statistic by using a different threshold
for each scale, which resulted in Corollary \ref{cor:allK}.

\begin{supplement}[id=supp]
\sname{Supplement}
\stitle{Technical Arguments}
\slink[doi]{10.1214/10-AOS839SUPP}
\sdatatype{.pdf}
\sfilename{cluster-suppl.pdf}
\sdescription{In the supplementary file \cite{clustersuppl}, we prove the results
stated here. It is divided into three sections. In the first section,
we state and prove general lower bounds on the minimax rate and upper
bounds on the detection rate achieved by an $\eps$-scan statistic. We
do this for the normal location model first and extend these results to
a general one-parameter exponential family. In the second section, we
gather a number of results on volumes and node counts. In the third and
last section, we prove the main results.}
\end{supplement}


%
\printaddresses

\end{document}